 \documentclass[12pt]{article}

\usepackage{blkarray}
\usepackage{epsfig,wrapfig,epic}
\usepackage{amsmath}
\usepackage{amscd}
\usepackage{wrapfig}
\usepackage{amsfonts}
\usepackage{mathrsfs}
\usepackage{color}

\usepackage{graphics}
\usepackage{graphicx}
\usepackage{amssymb}
\usepackage{amsmath}
\usepackage{url}

\DeclareMathAlphabet{\mathpzc}{OT1}{pzc}{m}{it}

\newcommand{\rank}{\mathrm{rank}}
\newcommand{\GIF}{\mbox{\texttt{gif}}}
\newcommand{\MMtx}{\mbox{\texttt{M}}}
\newcommand{\gen}{\mbox{\texttt{gen}}}

\newcommand{\Symbol}{\mathcal{S}}

\newcommand{\Pro}{\widehat{\textbf{D}}}
\newcommand{\ProRowsp}{\widetilde{\textbf{D}}}
\newcommand{\ProKer}{\textbf{D}}
\newcommand{\Proj}{\widehat{\boldsymbol{\pi}}}
\newcommand{\ProjRowsp}{\widetilde{\boldsymbol{\pi}}}
\newcommand{\ProjKer}{\boldsymbol{\pi}}

\newtheorem{algorithm}{Algorithm}[section]
\newtheorem{example}{Example}[section]

\newtheorem{theorem}{Theorem}[section]

\newtheorem{remark}[theorem]{Remark}
\newtheorem{define}[theorem]{Definition}

\title{
Geometric involutive bases for positive dimensional polynomial ideals and SDP methods  }

\author{ Greg
Reid and Fei Wang\thanks{Applied Mathematics Dept., Western
University, Canada. {\em Email:} reid@uwo.ca. Work partly funded by
Reid's NSERC grant (government of Canada).}
 \hspace{1cm} Wenyuan Wu\thanks{Chongqing Institute of Green and Intelligent Technology,
Chinese Academy of Sciences. {\em Email:} wuwenyuan@cigit.ac.cn.\;
The work is partly supported by the projects NSFC-11001040, West
Light Foundation of CAS from China (Y21Z010C10).}
}

\begin{document}

\input amssym.def
\maketitle

\begin{abstract}

Geometric involutive bases for polynomial systems of equations have their origin in the prolongation and projection methods of the geometers Cartan and Kuranishi for systems of PDE.
They are useful for numerical ideal membership testing and the solution of polynomial systems.
In this paper we further develop our symbolic-numeric methods for such bases.
We give methods to explicitly extract and decrease the degree
of intermediate systems and the output basis.
Algorithms for the numerical computation of involutivity criteria for positive dimensional
ideals are also discussed.

We were also motivated by some remarkable recent work by Lasserre and collaborators
who employed our prolongation projection involutive criteria as a part of their
semi-definite based programming (SDP) method for identifying the real radical of zero dimensional polynomial
ideals.  Consequently in this paper we begin an exploration of the interaction between geometric involutive bases
and these methods particularly in the positive dimensional case. Motivated by the extension of these methods to the positive dimensional case
we explore the interplay between geometric involutive bases and the new SDP methods.

\end{abstract}

\section{Introduction} 
\label{s:intro}

This paper is part of a stream devoted to developing symbolic-numeric prolongation projection algorithms for general systems  of partial and differential algebraic equations.
Such algorithms prolong (differentiate) such systems and project the prolonged systems to determine
obstructions or missing constraints to their integrability.
See Kuranishi \cite{Kuranishi:1957} for proof of termination of such methods using Cartan's
geometric involutivity criteria.
A by-product of these methods has been their implementation for linear
homogeneous partial differential equations with constant coefficients, and consequently for
polynomial algebraic systems.  See \cite{GerdtBlinkov:1998} for applications and symbolic algorithms for
polynomial systems.
The symbolic-numeric version of a geometric involutive form was first described and implemented in
Wittkopf and Reid \cite{ReidWittkopf:2001}.  It was applied to approximate symmetries of differential 
equations in \cite{BLRSZ:2004} and to polynomial solving in \cite{ReidZhi2009,ReidTangZhi:2003,SRWZ2010}.
See \cite{WuZhi12} where it is applied to the deflation of multiplicities in
multivariate polynomial solving.

The current paper is focused on further development of our geometric involutive basis algorithm
particularly in the positive dimensional case, and also in relation to real solving.
It is especially motivated by remarkable recent developments concerning real solution of
such systems by Lasserre, Laurent and Rostalski \cite{LasserreLaurentRostalski09} and their use of aspects of our prolongation projection algorithm
in the paper ``\textit{A prolongation-projection algorithm for computing the finite real variety of an ideal}''.
They developed a new approach for computing the real radical of zero dimensional polynomial
systems using semi-definite programming (SDP) techniques. 
See \cite{CurtoFialkow96}  for early fundamental work on such problems. 
Zero dimensional systems are those
having finitely many real solutions, and the real radical is the set of polynomials which vanish on these
solutions.  In contrast to the input systems the output radical systems from their approach are multiplicity free and
so are better conditioned for numerical solution techniques.   The output radical systems only have real roots and no complex roots.
This leads to possibility of lower complexity methods, since current methods for finding real solutions, mostly explicitly, or implicitly pass through complex root formulations.
Given the widespread popularity of linear programming (and by implication) SDP methods, the surprising
links between this area also open interesting research possibilities.
See \cite{BPT2012} for a recent book on the connections between semi-definite optimization and convex
algebraic geometry.

We briefly list some background references.
There have been considerable recent advances in numerical complex geometry.
See especially the books \cite{SommeseWampler05,BHSW13} and the references therein.
In approaches based on homotopy continuation, positive dimensional components characterize
the variety over $\mathbb{C}$ by certain witness points cut out by intersections of the components with random linear spaces.
For a modern text with many references on computational real algebraic geometry
see \cite{BasuPollackRoy06}.
Real algebraic geometry is a vast subject with many applications.
Sturm's ancient method on counting real roots of a polynomial in an
interval is central to Tarski's real quantifier elimination
\cite{Tarski31} and was further developed by Seidenberg
\cite{Seidenberg54}. One of the most important algorithms of real
algebraic geometry is cylindrical algebraic decomposition.  CAD was
introduced by Collins \cite{Collins75} and improved by Hong
\cite{Hong90} who made Tarski's quantifier
elimination algorithmic. This algorithm decomposes $\mathbb{R}^n$ into cells on
which each polynomial of a given system has constant sign.
The projections of two cells in $\mathbb{R}^n$ to
$\mathbb{R}^k$ with $k < n$ either don't intersect or are equal. The
computational cost of this algorithm, which is doubly exponential
\cite{DavenportHeintz88}, is a major barrier to its application.
See \cite{CMXY09} and \cite{Chen13} for modern improvements using
triangular decompositions.
For approaches based on obtaining witness points for the
real positive dimensional case see \cite{RRS00,GrigorevVorobjov92,Hauenstein12,WuReid13}.
Homotopy methods are used in \cite{Lu06}
and \cite{BDHSW12} for real algebraic geometry.
Recently such moment matrix completion techniques are explored by Zhi et al in \cite{MaZhi12} for finding at least one real root of a given semi-algebraic system.  Furthermore, based on critical point technique and moment matrix completion, they studied the computation of verified real solutions on positive dimensional system in \cite{YZZ13}.

As part of our initial exploration of this area, in this paper, we make some improvements
in our geometric involutive bases, by enabling the explicit extraction of projected systems
and hence reducing the size of matrices that can appear in intermediate computations.
Similarly motivated by the extension of these methods to the positive dimensional case
we explore the interplay between geometric involutive bases and the new SDP methods.
The symbol space of a polynomial system  or kernel of the matrix of its highest coefficients is the geometric generalization of the highest coefficient of
a polynomial.  Certain projections within the symbol space encode a geometric test - an analogue
of the S-polynomials in Gr\"obner basis approaches - for new members of the polynomial ideal.
We provide details and example of this in the numerical case.
An attempt in this paper is made to minimize use of terminology from the jet geometry of
partial differential equations, in order to make this accessible to a wider audience.

\section{Brief background on ideals and varieties}

\label{s:BriefBackGround}

In this section we briefly sketch some basic objects from real and complex algebraic geometry
and introduce some notation  for our paper.

\subsection{Some basic objects in complex algebraic geometry}
\label{s:CompAlgGeom}

Consider the set $\mathbb{C}[x_1, x_2, ... , x_n ]$
of multivariate polynomials with complex coefficients in the complex variables
$x = (x_1, x_2, ... , x_n) \in \mathbb{C}^n$.
Then $\mathbb{C}[x_1, x_2, ... , x_n ]$ is a ring.
Given $ P = \{p_1(x), p_2(x), ... , p_m(x) \} \subseteq
\mathbb{C}[x_1, x_2, ... , x_n ]=\mathbb{C}[x]$
its solution set or variety is:
\begin{equation}
\label{VC}
V_\mathbb{C} (p_1, p_2, ... , p_m) = \left\{ x \in \mathbb{C}^n :  p_j(x) = 0, 1 \leq j \leq m  \right\}
\end{equation}
For brevity we sometimes write $ V_\mathbb{C}  (P) = \{ x \in \mathbb{C}^n :  P(x) = 0  \}$.
Upper case letters $P$, $Q$, $R$, etc will denote sets of polynomials and lower case letters
$p$, $q$ etc will denote individual polynomials.

The ideal over $\mathbb{C}$ generated by $ P = \{ p_1,...,p_k \} $ is:
\begin{equation}\label{CIdeal}
    \left\langle P \right\rangle_\mathbb{C} = \left\langle p_1,...,p_k \right\rangle_\mathbb{C}  = \{ f_1 p_1 + ... + f_k p_k : f_j \in \mathbb{C}[x] ,    1\leq j \leq k\}
\end{equation}
and its associated radical ideal over $\mathbb{C}$ is
\begin{eqnarray}\label{RadCIdeal}
\sqrt[\mathbb{C}]{ \left\langle P \right\rangle }
		  &=&  \{ f   \in \mathbb{C}[x] : f(x) = 0
								\;  \mbox{ for all } \;  x \in V_\mathbb{C} (P) \}   \nonumber \\
              &=&
			 \{ f   \in \mathbb{C}[x] : f^m \in \left\langle P \right\rangle
								\;  \mbox{ for some } \;  m \in \mathbb{N} \}
\end{eqnarray}
where $\mathbb{N}$ is the set of non-negative integers.
\begin{example}\label{ex:simplestC}
To make this paper accessible to a wide audience we illustrate first some of the main
ideas on the simple and well-known case of systems of univariate polynomials.
Given a system of $k$ univariate polynomials  $P = \{ p_1 , ... , p_k \}$ with coefficients from some computable field (e.g. $\mathbb{Q}$), a Gr\"obner basis (or $\gcd$) computation
returns a single polynomial $q(x)$:\begin{equation}\label{SimplestCIdeal}
    \left\langle q \right\rangle_\mathbb{C} = \left\langle p_1,...,p_k \right\rangle_\mathbb{C}
\end{equation}
The factorization of $q(x)$ over $\mathbb{C}$ has form:
\begin{equation}\label{CpFac}
q(x) = a (x - a_1)^{n_1} ... (x - a_\ell )^{n_\ell }
\end{equation}
where the roots $a_j \in \mathbb{C}$ of $q(x)$ are distinct.
Though the $a_j$ can't be found in general by finitely many rational operations
the so-called square-free factorization can be found by such operations yielding:
\begin{equation}\label{CpFac}
\tilde{q}(x) =  \frac{q(x)}{\gcd (q(x), q'(x) )} =  a (x - a_1) ... (x - a_\ell )
\end{equation}
For this example the ideal, variety and radical ideal over $\mathbb{C}$ are:
\begin{eqnarray}\label{SimpleIdealCEx}
  \left\langle P \right\rangle_\mathbb{C} &=& \{ g(x) \cdot   (x - a_1)^{n_1} ... (x - a_\ell )^{n_\ell }  : g(x) \in \mathbb{C}[x]  \}    \nonumber \\
                V_\mathbb{C} (P) &=& \{ a_1 , a_2 , ... , a_\ell \} \\
\sqrt[\mathbb{C}]{ \left\langle P \right\rangle }
		  &=&  \{  g(x)  \cdot (x - a_1) ... (x - a_\ell ): g(x) \in \mathbb{C}[x]   \}   \nonumber
\end{eqnarray}
For sophisticated generalizations to primary decomposition for multivariate systems see Gianni et al. \cite{GTZ88}.

\end{example}

\subsection{Some basic objects in real algebraic geometry}
\label{s:RealAlgGeom}

Suppose that $x = (x_1, x_2, ... , x_n ) \in \mathbb{R}^n$
and consider a system of $k$ multivariate polynomials
$P = \{ p_1(x), p_2(x), ... , p_k(x) \}  \subseteq
\mathbb{R}[x_1, x_2, ... , x_n ]$ with
real coefficients.
Its solution set or variety is
\begin{equation}\label{VReal}
    V_\mathbb{R}  (p_1,...,p_k) = \{ x \in \mathbb{R}^n :  p_j(x) = 0,\; 1\leq j \leq k\}
\end{equation}

The ideal generated by $ P  =  \{   p_1,...,p_k  \}  \subseteq \mathbb{R}$  is:
\begin{equation}\label{RIdeal}
     \left\langle P  \right\rangle_\mathbb{R} = \left\langle  p_1,...,p_k  \right\rangle_\mathbb{R}  = \{ f_1 p_1 + ... + f_k p_k : f_j \in \mathbb{R}[x] ,    1\leq j \leq k  \}
\end{equation}
and its associated radical ideal over $\mathbb{R}$ is defined as
\begin{eqnarray}\label{RadRIdeal}
\sqrt[\mathbb{R}]{ \left\langle P \right\rangle }
              &=&
			 \{ f   \in \mathbb{R}[x] : f^{2m}+\Sigma_{j=1}^{s} q_j^{2} \in \left\langle P \right\rangle
								\;  \mbox{for some} \; q_j  \in \mathbb{R}[x], m \in \mathbb{N} \backslash \{0\} \} \label{defrrad}
\end{eqnarray}
A fundmental result \cite{BasuPollackRoy06} is:
\begin{theorem} \label{RealNull}
[Real Nullstellensatz]   For any ideal $I \subseteq \mathbb{R}[x]$ we have $\sqrt[\mathbb{R}]{ I } = I(V_R(I))$.
\end{theorem}
 Consequently
\begin{eqnarray}\label{RadRIdeal}
\sqrt[\mathbb{R}]{ \left\langle P \right\rangle }
              &=&  \{ f(x)   \in \mathbb{R}[x] : f(x) = 0
								\;  \mbox{ for all } \;  x \in V_\mathbb{R} (P) \} \;
\end{eqnarray}

 \begin{remark} \label{radical}
 An ideal $I \subseteq \mathbb{R}[x]$ is real radical if and only if
     for all $ p_1, \cdots, p_k \in R[x]$:
 \begin{equation}
  p_1^2 + \cdots + p_k^2 \in I \Longrightarrow p_1, \cdots, p_k \in I.
 \end{equation}
 \end{remark}
For these and many other results see  \cite{BasuPollackRoy06} and the references
cited therein.

 \begin{example}\label{ex:simplestR}
Consider the simplest case of a system of $k$ univariate polynomials in  some computable subfield
of  $\mathbb{R}$ (e.g. $\mathbb{Q}$).   Then as in the complex case a Gr\"obner basis
of such a system yields a single polynomial $q(x)$ having the same roots.
Discarding the factors with complex roots with nonzero imaginary parts yields a polynomial of form:
\begin{equation}\label{RpFac}
\tilde{q}(x) =  b (x - b_1)^{m_1}  ... (x - b_j )^{m_j}
\end{equation}
where $b_1$, $b_2$, ... , $b_j$ are the real roots and $m_1$, ... , $m_j$ their corresponding multiplicities.
Then
\begin{eqnarray}\label{SimpleIdealREx}
  \left\langle P \right\rangle_\mathbb{R} &=& \{ f(x) \cdot   (x - b_1)^{m_1} ... (x - b_j )^{m_j }  : f(x) \in \mathbb{R}[x]  \}  \nonumber  \\
                V_\mathbb{R} (P) &=& \{ b_1 , b_2 , ... , b_j \} \\
\sqrt[\mathbb{R}]{ \left\langle P \right\rangle }
		  &=&  \{  g(x) \cdot  (x - b_1) ... (x - b_j ) : g(x) \in \mathbb{R}[x]   \}   \nonumber
\end{eqnarray}

\end{example}

\section{Geometric prolongation and projection for polynomial systems}
\label{s:GeoProProj}

In this section we give a brief description the well-known presentation of
polynomial systems as linear functions of their monomials and the related
coefficient matrix and its kernel and rowspace \cite{Stetter:2004,Mourrain:1996,Mourrain:1999,MollerSauer:2000} and historical work by Macaulay \cite{Macaulay:1916}.
We describe a type of elimination called geometric projection
and then describe geometric prolongation resulting from multiplying polynomials
by monomials.

We exploit the well-known correspondence between polynomial
systems and systems of constant coefficient linear homogeneous
PDE. This equivalence has been extensively studied and exploited
in the exact case by Gerdt \cite{GerdtBlinkov:1998} and his
co-workers in their development of {\it involutive bases}.
Our geometric involutive bases are involutive by the geometric criteria in
\cite{Kuranishi:1957,Pommaret:1978,Seiler2010} and more distantly related
to that of \cite{GerdtBlinkov:1998} which are closer relatives of Gr\"obner bases.

Consider a system of $\ell$ polynomials $P \subseteq \mathbb{K}[x]$ of degree $d$
in the variables $x = (x_1, ... , x_n )$ where $\mathbb{K} = \mathbb{R}$ or $\mathbb{C}$.
Monomials are denoted by $x^\alpha := x_1^{\alpha_1} ...  x_n^{\alpha_n}$ where $\alpha \in \mathbb{N}$
and  the degree of $x^\alpha$ is $|  \alpha| = \alpha_1 + ... + \alpha_n$.
Then the system $P$ can be written as:
\begin{equation}
\label{sigPoly}
P  =  \left\{ \sum_{ | \alpha|  \leq d} \;   a_{k, \alpha} \: x^\alpha :  k = 1 ... \ell   \right\}
\end{equation}
To apply the methods of numerical linear algebra the system is converted into matrix form
\cite{Stetter:2004,Mourrain:1996,Mourrain:1999,MollerSauer:2000}.

\begin{define}[Coefficient Matrix $\mathbf{C(P)}$, $\mathbf{J^d}$ and vector of monomials]
\label{defCoeffMtx}
Denote the coefficient matrix
of $P$ in (\ref{sigPoly}) by $C(P)$.  Let $\textbf{x}^{(\leq d )}$ be the column vector of monomials $x^\alpha$
with $0 \leq \alpha \leq d$ sorted by total degree.  We suppose that the columns of
$C(P)$ are sorted in the same order.  Then $P = C(P) \textbf{x}^{(\leq d )}$
where $C(P) \in \mathbb{R}^{\ell \times N(n,d) }$ and
 $N(n,d) :=
 \small{\left( \begin{array}{c}
    d + n \\
       d \\
 \end{array} \right) } $ is the number of monomials in $\textbf{x}^{(\leq d )}$.
Polynomials can be equivalently represented by the row vectors of $C(P)$, that is as vectors
in $J^d := \mathbb{R}^{N(n,d)}$.
\end{define}
Prolonging polynomials by multiplying them by monomials is an essential geometric operation in this
paper.
\begin{define}[prolongations  $\Pro$ and $\ProRowsp$]
Consider a system of polynomials $P$ of degree $d$.
Let $p \in P$ have degree $\bar{d}$.  Then the prolongation of
$p$ written $\Pro(p)$ is defined as $\Pro (p) = \{ p \} \cup \{ x_j p : 1 \leq j \leq n \}$.
The prolongation of the system $P$ is defined as
$\Pro^k  (P) =  \{  x^\alpha p :  0 \leq \deg (x^\alpha p ) \leq d+k,  \alpha \in \mathbb{N}^n, p \in P \}$.
Equivalently we can represent the prolongation geometrically as the span of
the corresponding row vectors of $C(\Pro^k P )$, which we denote by
$\ProRowsp^k (P)  := \mbox{rowsp} (C(\Pro^k P ) )$ which is a subspace
of $J^{d+k}$.
\end{define}

\begin{example}
Suppose $x = (y, z)$ and $P = \{  2,  2 y + z \}$.
 Then $\Pro (P) = \{  2, 2y , 2z, 2y^2, 2 yz, 2z^2,  \\
 2y^2 + yz, 2yz + z^2 \}$.
\end{example}

\begin{define}[projections $\Proj$ and $\ProjRowsp$]
Consider a polynomial system of degree $d \geq 1$ written in the form
$P = C(P) \textbf{x}^{(\leq d )}$ with the columns of $C(P)$ sorted
in descending order by degree.
 The rows in the Gauss echelon form of $C(P)$ with pivots of degree
less than $d$ span a subspace of $J^{d-1}$ which we denote
by $\ProjRowsp (P)$.
We denote the set of polynomials of degree $ \leq d-1$ corresponding
to the row vectors by $\Proj (P)$.
Iterations of projections $\Proj^\ell (P) \subset \mathbb{R}[x]$
 and equivalently $\ProjRowsp^\ell (P) \subset J^{d - \ell} $ are defined similarly.
\end{define}
We have adopted an abbreviated notation for prolongation and projection here to avoid cumbersome
indices indicating the spaces on which these operators act.
We will also need to prolong and project kernels of the coefficient matrices
of polynomial systems.
\begin{define}[prolongation $\ProKer$ and projection $\ProjKer$ on the kernel]
Consider a polynomial system $P \subset \mathbb{R}[x]$ of degree $d$.
Given a subspace $V$ of $J^d$ and $\ell \leq d$
define $\ProjKer^\ell (V)$ as the vectors of $V$ with
the components of degree $\geq d - \ell$ discarded.
To abbreviate notation we will write
 $\ProjKer^\ell (P) := \ProjKer^\ell  \ker C(P)$.
The $k$-th prolongation of the kernel is
$\ProKer^k  (P) :=  \ker C(\Pro^k P)$.
\end{define}
In summary we have presented three (!) notations for prolongation and projection since we need
to work directly with them sometimes as polynomial systems,  and sometimes row spaces or kernels.
The row space and kernel are orthogonal to each other in $J^d$.
Projection in the kernel is the usual projection operator $\ProjKer^\ell$.
Geometrically  the corresponding projection in the row space can be obtained as the orthogonal
complement of  $\ProjKer^\ell (P)$.  Alternatively it can be obtained by first considering
$J^{d - \ell}$ as a subspace in $J^d$ and then intersecting the subspace $J^{d - \ell}$ with
$\mbox{rowsp} (P)$.

Suppose that $A = C(P)$ is the coefficient matrix of a system of polynomials $P$.
To numerically implement an approximate involutive form method, we
proposed  in \cite{ BLRSZ:2004,ReidWittkopf:2001,ReidTangZhi:2003} a numeric version of the projection
operator based on singular value decomposition (SVD).
We first find the SVD of $A$ given by $A = U \cdot \Sigma \cdot V $
where $U$ and $V$ are unitary matrices and $\Sigma$ is a diagonal matrix
whose diagonal entries are real decreasing non-negative numbers.
The approximate rank $r$ is the number of singular values bigger
than a fixed tolerance.  Deleting the first $r$ rows of $V$ yields
an approximate basis for $\ker A$ and
an estimate for $\dim \ker A$. Deleting highest degree components of the vectors in this basis, yields an approximate spanning set for $\ProjKer \ker A$  and an estimate for $\dim \ProjKer \ker A$. If desired further computation yields bases for $\ProjKer \ker A$. Then we compute the kernel of the spanning set
of $\ProjKer \ker A$.  Similarly we can compute approximate spanning sets and if desired bases of the prolongations and projections of the system.

 \begin{remark}[Alternative representations and extraction of intermediate systems]
 \label{BijectionGenKer}
\
\\
In summary
 prolongation and projection can equivalently be computed in either the kernel, the rowspace, and at any time polynomial
generators can be extracted. Underlying this is a 1 to 1 correspondence between vector spaces (not elememts): in particular between the row spaces and its orthogonal complement, the kernel. 

\end{remark}

\begin{example}
\label{ex:Deg8Univar}
Consider
\begin{equation}
\label{eq:DeguUnivar}
P = \{ x^8  -  x^4 - 2 ,   x^8  - 3 x^4  + 2 \}  \subseteq \mathbb{R}[x]
\end{equation}
Here the coefficient matrix is given by $C(P)$ below:
\begin{equation}
\label{CMtxSimpleEx}
C(P) \cdot \mbox{\textbf{x}}^{(\leq 8 )}=    \left( \begin{array}{ccccccccr}
     1 & 0 & 0 & 0 & -1 & 0 & 0 & 0 & -2 \\
     1 & 0 & 0 & 0 & -3 & 0 & 0 & 0 &  2 \\
\end{array} \right)
\small{ \left[ \begin{array}{c}
	x^8 \\
	x^7 \\
	\vdots \\
	x^1 \\
	 1  \\
\end{array} \right]
 }
=
\small{
 \left( \begin{array}{c}
        0   \\
	 0  \\
\end{array} \right)
}
\end{equation}
The most familiar computation for most readers is to eliminate the polynomials as in
a Gr\"obner basis calculation:  $ x^8  -  x^4 - 2 - (x^8  - 3 x^4  + 2) =
2 x^4 - 4$.  This can also be done as a computation on the row space of
$C(P)$, yielding the result as the generator of $\ProjRowsp^4  P$.
Equivalently by Remark \ref{BijectionGenKer} we can compute the result by
projecting basis vectors of the kernel of $C(P)$ obtaining
$\ProjKer^4 P $ and then recover the generator $2 x^4 - 4$.
The original $8$ degree polynomials can be discarded since they are consequences
of $2 x^4 - 4$.  In particular a Gr\"obner basis for the ideal generated by $P$ is
\begin{equation}
\label{GBsimpex}
 x^4 - 2
\end{equation}

The kernel of $C(P)$ is easily calculated numerically by the SVD.
We obtain the table of dimensions for the projections of $\ker C(P)$
in Figure \ref{TableCMtxSimpleEx}.
We use singular value decomposition to compute its kernel and then project its
vectors to $\ProjKer^4 P$.  The generator corresponding to this projection is:
\begin{equation}
 0.4472136\,{x}^{4}- 0.8944272.
\end{equation}
where the coefficients here and elsewhere in the paper have been truncated from 15 digits to 7 digits.
After normalization, we get the generator $x^4 - 2$.

\begin{figure}
\centerline{
\begin{tabular}{|l|cccc|}
\hline
            & $k=0$ & $k=1$ & $k=2$ & k=3    \\
\hline
$\ell=0$ &  7       &  6      &  5        &    4     \\
$\ell=1$ &  7       &  6      &  5        &    4    \\
$\ell=2$ &  6       &  6      &  5        &    4    \\
$\ell=3$ &  5       &  5      &  5        &   \color{blue}{\boxed{4}}    \\
$\ell=4$ & \color{red}{\boxed{4}}      &  4      &  4        &    4    \\
$\ell=5$ &  4       &  4      &  4        &    4    \\
$\ell=6$ &  3       &  4      &  4        &    4    \\
$\ell=7$ &  2       &  3      &  4        &    \boxed{4 }   \\
$\ell=8$ &  1       &  2      &  3        &    4    \\
$\ell=9$ &           &  1      &  2        &   3     \\
$\ell=10$&          &          &  1        &   2     \\
$\ell=11$&          &          &            &   1     \\
\hline
\end{tabular}}
\caption{ Table of $\dim   \ProjKer^\ell   \ProKer^k P$ for
        (\ref{CMtxSimpleEx}) for Example \ref{ex:Deg8Univar}.  The (red) boxed $4$ in the first
		column corresponds to $\ProjKer^4 P$ and a geometric involutive
		basis for $P$ as found by Algorithm \ref{alg:ProjInvBasis}. 
		The blue and black boxed $4$'s in the fourth
		column correspond to geometric involutive
		bases for $P$.
		\label{TableCMtxSimpleEx}}
\end{figure}

\end{example}

\section{Geometric involutive bases}
\label{s:GeoInvForm}

In this section we describe geometric involutive form.
For a more detailed description see
\cite{BLRSZ:2004, ReidSmithVerschelde:2002,ReidWittkopf:2001,
ReidTangZhi:2003}.

Exact elimination methods for exactly given polynomial systems
(e.g. Gr\"obner Bases), usually employ Gaussian Elimination (e.g.
linear elimination of monomials). Such exact methods usually
depend on the ordering of input (e.g. term ordering in the case of
Gr\"obner Bases), and so are coordinate dependent. Since the order
of elimination can force division by small leading entries, such
methods are generally unstable, when used on approximate systems.
In contrast, exact elimination methods from the geometric theory
of PDE are coordinate independent \cite{Kuranishi:1957,Pommaret:1978}
and this motivated our study of numerical versions of such methods
which is continued in this paper.

\subsection{Symbol, class and Cartan involution test}

\begin{define}[Symbol matrix and class of a monomial]
Given a polynomial system of degree $d$, its symbol matrix, denoted $\Symbol (P)$ is the submatrix
of $C(P)$ corresponding to its degree $d$ monomials.
Consider a monomial $x^\alpha$ where $\alpha = (\alpha_1, ... , \alpha_n ) \in \mathbb{N}^n$.
Then the class of $x^\alpha$ is the least $j$ such that $\alpha_j \not = 0$.
\end{define}
For Example \ref{ex:Deg8Univar} the symbol matrix is
the submatrix
 $\left( \begin{array}{c}
     1  \\
     1  \\
\end{array} \right)  $
 of $C(P)$ given in (\ref{CMtxSimpleEx}).
Consider the system
\begin{equation}
\label{simpSymSys}
P = \{ {x_2}^2 - 1 , 2  x_1 x_2 - 3 x_1  \}
\end{equation}
For what follows we sort the columns of the symbol matrix in descending order according to class.
The degree two monomials are $x^{(0,2)} = x_2^2$, $x^{(1,1)} = x_1 x_2 $,  $x^{(2,0)} = x_1^2$.  Here $x_2^2$ is class $2$.
Monomials $x_1 x_2$ and $x_1^2$ are class $1$.
Then the symbol matrix is:
\begin{equation}
\label{SymSimpMtx}
\Symbol (P) = \left( \begin{array}{ccc}
     1 &  0 &  0  \\
     0 &  2 & 0 \\
\end{array} \right)
\end{equation}

\begin{define}[Cartan test for involutivity of the Symbol]
Suppose that the columns of the symbol matrix for a system of degree $d$ are sorted in descending order by class and that
it is reduced to Gauss echelon form.  For $k = 1, 2, ... , n$ define the quantities $\beta_k$ as the number of pivots in this reduced matrix of class $k$.  Then in a generic system of coordinates
the symbol is involutive if:
\begin{equation}
\label{CartanTest}
			\sum_{k=1}^{k=n} k \beta_d^{(k)} = \rank \: \Symbol ( \Pro P )
\end{equation}

\end{define}

The following combinatorial quantities will be useful in our numerical determination
of involutivity of symbol matrices.
Consider systems in $n$ variables of degree $d$.

Denote:
\begin{equation}
\begin{array}{lcccccc}
N(n,d) & = &   \left( \begin{array}{c}
							    n + d \\
							       d \\
 				\end{array} \right)  & =& \mbox{Number of monomials of degree} \leq  d  \\
N_{\mbox{deg}}  (n,d) & = &   \left( \begin{array}{c}
							    n + d - 1 \\
							       d \\
 					\end{array} \right)  & =& \mbox{Number of monomials of degree} \;  d   \\
N_c (n,d,k) & = &   \left( \begin{array}{c}
							    n + d - k - 1 \\
							          d - 1 \\
 					\end{array} \right)  & =& \mbox{Number of class $k$ monomials of degree} \; d 
\end{array}
\end{equation}

\begin{example}
\label{involSimpEx}
For system $P$ given in (\ref{simpSymSys}):
\begin{equation}
\label{simpnums}
N(2,2) = 6,  N_{\mbox{deg}} (2,2) = 2,  N_c(2,2,1) = 2,  N_c(2,2,2) = 1
\end{equation}
The symbol matrix (\ref{SymSimpMtx})  is already in Gauss echelon form with respect to class.  There is one pivot of class $2$ so $\beta_2^{(2)} = 1$ and one pivot of class $1$ so
$\beta^{(1)}_2 = 1$.  Also an easy calculation gives $ \rank \:  \Symbol ( \Pro P )  = 3$.
So
\begin{equation}
\label{CartanTestSimpEx}
			\sum_{k=1}^{k=2} k \beta_d^{(k)} = 3 = \rank  \:  \Symbol ( \Pro P )
\end{equation}
and the symbol is involutive.  In all cases $\sum_{k=1}^{k=n} k \beta_d^{(k)} \leq \rank \:  \Symbol ( \Pro P ) $.  Indeed in our example if we reverse the order of the coordinates and
recalculate we get
$\Symbol (P) =  \small{ \left( \begin{array}{ccc}
     0 &  2 &  0  \\
     0 &  0 & 1 \\
\end{array} \right) }$.
Then $\beta_2^{(2)} = 0$, $\beta^{(1)}_2 = 2$ and
$\sum_{k=1}^{k=2} k \beta_d^{(k)} = 2 < \rank \:  \Symbol ( \Pro P ) $
so the test indicates a non-involutive symbol however the result may be due
to the coordinates being nongeneric which is indeed the case here.
 A generic linear change of coordinates
by a random $2 \times 2$ matrix then shows the symbol is involutive.
\end{example}

To extract a matrix for the symbol space of the variables of degree $d$ we proceed as follows
for a system $P$ of degree $d' \geq d$.
Suppose that vectors that are a basis for the kernel of $C(P)$ form the rows of a matrix $B$.
First numerically
project the kernel of the system $P$ onto the subspace $J^d$ via $\ProjKer^{d' - d} P$
by deleting the coordinates in the basis of degree $> d$ to obtain for  $\ProjKer^{d' - d} P$ a spanning set $\tilde{B}$ given by the remaining rows of $B$.
Then delete the columns in $\tilde{B}$
corresponding to variables of degree $< d$ to obtain a matrix $A_d$ corresponding to the
orthogonal complement of the symbol
for degree $d$.
Let $A_d^{(k)}$ be the submatrix of $\tilde{B}$ with columns corresponding to class $k$ or less deleted.
In generic coordinates
\begin{equation}
\label{beta[k]}
\beta_d^{(k)}= N_c (n, d, k) -  \left( \rank \:  A_{d}^{(k-1)} - \rank \: A_d^{(k)} \right) , \; \; \; k = 1 \ldots  n .
\end{equation}
Then the SVD can approximate the ranks in this equation for carrying out the Cartan Test (\ref{CartanTest}).

\begin{define}
[Involutive System] A system of polynomials $P \in \mathbb{R}[x]$ is involutive if
 $ \dim \:   \ProjKer \ProKer P = \dim \: P  $
and the symbol of $P$ is involutive.
\end{define}

\begin{define}[Projected Involutive System]
\label{def:projinvolsys}
Consider a system of polynomials $P \in \mathbb{R}[x]$ of degree $d$.
Suppose that $k$, $\ell$ are integers with $k \geq 0$ and $0 \leq \ell \leq k +d$.
 Then $\ProjKer^\ell \ProKer^k P$ is
projectively involutive at prolongation order $k$ and projected order $\ell$,
if $\ProjKer^\ell  \ProKer^k P$ satisfies
the projected elimination test
\begin{equation}
\label{projelimtest} {\dim} \; \ProjKer^\ell \ProKer^k P  = {\dim}\;
\ProjKer^{\ell+1} \ProKer^{k+1} P
\end{equation}
and the { symbol} of $\ProjKer^\ell \ProKer^k P $ is involutive.
\end{define}
In \cite{BLRSZ:2004} we prove:
\begin{theorem}
\label{projinvol<=>invol}
A system is projectively involutive if and only if it is involutive.
\end{theorem}
\begin{theorem}[Criterion for zero dimensional involutive system]
\label{secondtheorem}
A zero dimensional system of polynomials $P \in \mathbb{R}[x]$ is
 projectively involutive at order $k$ and projected order $\ell$
 if and only if $ \ProjKer^\ell  \ProKer^k P$ satisfies the projected
elimination test (\ref{projelimtest}) and
\begin{equation}
\label{projsymboltest} {\dim} \; \ProjKer^\ell \ProKer^k  P  = {\dim}
\;   \ProjKer^{\ell+1}  \ProKer^k P
\end{equation}
\end{theorem}
This criterion is used by Lasserre et al \cite{SemiDefinite:2012} in their prolongation projection
algorithm to determine the finite real radical.
When there are $2$ variables then it is easily shown that:
 \begin{equation} \label{projelimtest1}
 \Symbol \: \ProjKer^\ell \ProKer^k P  \; \;
\mbox{is involutive} \Longleftrightarrow {\dim} \: \Symbol
\:  \ProjKer^\ell  \ProKer^{k+1}  P =
{\dim} \: \Symbol \: \:
\ProjKer^{\ell}  \ProKer^k P
\end{equation}
 and this gives a computationally easy characterization by using
\begin{equation}
\label{symboleqn}{\dim} \:  {\Symbol} \:
\ProjKer^\ell   \ProKer^k  P ={\dim}
\:   \ProjKer^\ell  \ProKer^k  P  -{\dim} \:
\ProjKer^{\ell+1}  \ProKer^k P
\end{equation}
The criteria in (\ref{projelimtest}) applies to both zero and positive dimensional bivariate systems.
\subsection{Projected involutive form algorithm} \label{Sec:proinv}

The following method completes systems to approximate
involutive form.
We seek the smallest $k$ such that there exists an $\ell$
  with $ \ProjKer^\ell \ProKer^k P$ approximately involutive, and generates the same ideal as the input system. We choose
the system corresponding to the largest such $\ell  \leq  k$ if there are several such values for the
given $k$.

\begin{center}

\begin{algorithm}[Projected involutive basis]
\label{alg:ProjInvBasis}

\noindent

\begin{itemize}
\item[]    \parskip0mm  \textbf{Input:}
			   $ Q \subseteq \mathbb{R}[x_1,\ldots,x_n]$. A tolerance $\epsilon$.
			\item[] Set $k := 0$, $ d := \deg(Q)$ and $P := \ker C(Q)$
			\item[]  \textbf{repeat}
		           \begin{enumerate}					
				\item[] Compute $\ProKer^k (P)$
				\item[] Initialize set of involutive systems $I:=  \{ \}$					
                	\item[] \textbf{for} $\ell = 0 \cdots  (d+k) $ \textbf{do}
					\begin{itemize}
							\item[]  Compute $R := \ProjKer^\ell \ProKer^k (P)$
                           		\item[]\textbf{ if } $R$ involutive  \textbf{then} $I := I \cup \{ R \} $  \textbf{end if}
					\end{itemize}
					\item[]   \textbf{end do}
					\item[] Remove systems $\bar{R}$ from $I$ not satisfying $\ProKer^{d + k-\bar{d}} \bar{R}  \subseteq \ProKer^k (P)$
					\item[] $k := k + 1$
                     \end{enumerate}
			\item[] \textbf{until} $I \not = \{ \}$

\item[] {\bf Output}:  Return the polynomial generators  of the  involutive system $\bar{R}$ in $I$
\item[]  \hspace{2cm} of lowest degree $\bar{d}$.

\end{itemize}

\end{algorithm}

\end{center}

Note that this algorithm works on kernels, but could by Remark \ref{BijectionGenKer} equivalently
work on their orthogonal complements -- the associated row spaces.
The condition $\ProKer^{d + k-\bar{d}} \bar{R}  \subseteq \ProKer^k (P)$ is a standard subspace inclusion
test for the prolonged kernels.  It ensures that the output system generates the same ideal as the
input system and has the same solutions. 

\begin{remark}[Decreasing degrees by extracting involutive projections]
\label{rem:extracProj}  \
\\ We note that a simple illustration of Algorithm \ref{alg:ProjInvBasis} is Example \ref{ex:Deg8Univar} where all univariate polynomials are involutive. This algorithm is an improvement on that published in \cite{SRWZ2010} where to ensure the inclusion conditions for positive dimensional ideals, the number of projections was limited to $0 \leq \ell \leq k$.  
\end{remark}

\section{Moment matrices and SDP }

\label{s:MM-SDP}

\subsection{Moment Matrices}

\label{s:MM}

Here we focus just on the construction of moment matrices.  For the theoretical background the reader is directed to \cite{SemiDefinite:2012}.

A moment matrix is a symmetric matrix $M = (M_{\alpha, \beta})$ indexed by $\mathbb{N}^{n}$ ($\alpha, \beta \in\mathbb{N}^n $).  Here $\alpha$ is the index for rows, $\beta$ is the index for columns.    Without loss $M_{0,0} = 1$.

Given a multivariate polynomial system $P \subseteq \mathbb{R}[x_1, ... , x_n]$.
Let  $d = \deg(P)$ and $M  \in \mathbb{R}^{N(n,d) \times N(n,d)}$ be the truncated moment matrix.
The linear constraints imposed by $P$ are constructed as
\begin{eqnarray} \label{MMtxLNC}
M \cdot A^{T} =0; \; A = C (\Pro^d (P)),
\end{eqnarray}
where $C$ is the coefficient matrix function given in Definition \ref{defCoeffMtx}.

\subsection{Moment matrix for univariate example}

\label{ex:simple}

In Example \ref{ex:Deg8Univar} a degree $8$ input system was reduced to a degree $4$
output polynomial $p = x^4 - 2$.
Then in matrix form the polynomial is
\begin{equation}\label{simplesys}
    B v = \left(
            \begin{array}{ccccc}
              -2 & 0 & 0 & 0 & 1 \\
            \end{array}
          \right)
          \left(
            \begin{array}{c}
              u_0 \\
              u_1 \\
              u_2 \\
             u_3 \\
             u_4 \\
            \end{array}
          \right)
          = 0,
          \ker B = \mbox{span}_\mathbb{R}
		 {\tiny{
           \left\{  \left( \begin {array}{c} 1\\ \noalign{\medskip}0
\\ \noalign{\medskip}0\\ \noalign{\medskip}0\\ \noalign{\medskip}2
\end {array} \right) , \left( \begin {array}{c} 0\\ \noalign{\medskip}0
\\ \noalign{\medskip}0\\ \noalign{\medskip}1\\ \noalign{\medskip}0
\end {array} \right) , \left( \begin {array}{c} 0\\ \noalign{\medskip}0
\\ \noalign{\medskip}1\\ \noalign{\medskip}0\\ \noalign{\medskip}0
\end {array} \right) , \left( \begin {array}{c} 0\\ \noalign{\medskip}
1\\ \noalign{\medskip}0\\ \noalign{\medskip}0\\ \noalign{\medskip}0
\end {array} \right)  \right\}
		 }}
\end{equation}
		
The moment matrix is the infinite matrix whose $(\alpha, \beta)$ entry
is $u_{\alpha + \beta}$ and $\alpha$, $\beta$ $\in \mathbb{N}^n$
given by:

\begin{equation}\label{MMtx4}
M =
\left(
  \begin{array}{cccccc}
    u_0 & u_1 & u_2 & u_3 & u_4 & \cdots  \\
    u_1 & u_2 & u_3 & u_4 & u_5 & \cdots \\
    u_2 & u_3 & u_4 & u_5 & u_6 & \cdots \\
    u_3 & u_4 & u_5 & u_6 & u_7 & \cdots  \\
    u_4 & u_5 & u_6 & u_7 & u_8 & \cdots \\
    \vdots & \vdots & \vdots & \vdots & \vdots & \ddots \\
  \end{array}
\right)
\end{equation}
In the SDP-moment matrix approach the given polynomial system, in this case
$\{ x^4 - 2 \}$, is first prolonged to twice its degree:
\begin{equation}
\label{ProSimp8}
   \Pro^4 \{ x^4 - 2 \} =  \{ x^4 - 2, x^5 - 2x,  x^6 - 2x^2, x^7 - 2x^3, x^8 - 2 x^4 \}
\end{equation}
From (\ref{MMtxLNC})  the constraint system
when we impose $u_0 = 1$ is equivalent to the linear system
\begin{equation}
u_4 - 2 = 0,
u_5 - 2 u_1 = 0,
u_6 - 2 u_2 = 0,
u_7 - 2 u_3 = 0,
u_8 -  2 u_4  = 0
\end{equation}
which can be regarded as the rewrite rules:
$u_4 \rightarrow 2,
u_5 \rightarrow 2 u_1 ,
u_6 \rightarrow 2 u_2,
u_7 \rightarrow 2 u_3 ,
u_8 \rightarrow 2 u_4 \rightarrow 4$.
Imposing these constraints the truncated moment matrix to degree $8$ is
\begin{equation}\label{MMtx4c}
M =
\left(
  \begin{array}{ccccc}
     1   & u_1 & u_2 & u_3 & 2    \\
    u_1 & u_2 & u_3 & 2    & 2u_1 \\
    u_2 & u_3 & 2    & 2u_1 &2u_2 \\
    u_3 & 2    & 2u_1 & 2u_2 & 2u_3 \\
    2    & 2u_1 & 2u_2 & 2u_3 & 4   \\
  \end{array}
\right)
\end{equation}

The moment matrix~\eqref{MMtx4c} is then sent to the SDP solver Yalmip in Matlab
to numerically compute a generic point $(u_1, u_2, u_3)$ if possible such that
$M$ is a positive semidefinite matrix with maximum rank.
This solver returns an approximation which can be recognized for illustrative convenience
as $(u_1, u_2, u_3) = (0, \sqrt{2}, 0)$.
Its associated moment matrix and moment matrix kernel are:
\begin{equation}
	M =  \left( \begin {array}{ccccc} 1&0&\sqrt {2}&0&2\\ \noalign{\medskip}0&
\sqrt {2}&0&2&0\\ \noalign{\medskip}\sqrt {2}&0&2&0&2\,\sqrt {2}
\\ \noalign{\medskip}0&2&0&2\,\sqrt {2}&0\\ \noalign{\medskip}2&0&2\,
\sqrt {2}&0&4\end {array} \right) ,
\ker M =  \mbox{span}_\mathbb{R}\left\{   \left( \begin {array}{c}
-2\\ \noalign{\medskip}0\\ \noalign{\medskip}0\\ \noalign{\medskip}0
\\ \noalign{\medskip}1\end {array} \right) ,
\left( \begin {array}{c} -\sqrt {2}
\\ \noalign{\medskip}0\\ \noalign{\medskip}1\\ \noalign{\medskip}0
\\ \noalign{\medskip}0\end {array} \right) , \left( \begin {array}{c} 0\\ \noalign{\medskip}-\sqrt {2}
\\ \noalign{\medskip}0\\ \noalign{\medskip}1\\ \noalign{\medskip}0
\end {array} \right)  \right\}
\end{equation}
The kernel corresponds to the generating set
\begin{equation} \label{sdpker}
 \{ \sqrt{2} - x^2,  2 - x^4,  \sqrt{2}x - x^3 \}
\end{equation}
Applying geometric involutive form algorithm yields a geometric involutive basis 
\begin{equation} \label{genset}
\{  \sqrt{2} - x^2  \}
\end{equation}
The last two polynomials are a consequence of $ \sqrt{2} - x^2 $ by our inclusion test, so are discarded.
 By Rostalski \cite{SemiDefinite:2012}, this is a basis of the real radical.


\section{Combining geometric involutive bases and moment matrix methods}

\label{sec:combnGIF-MMtx}

\subsection{Geometric involutive form and moment matrix algorithms}
\label{sec:algorithmsGIF-MMtx}

In this section we outline algorithms for combining geometric involutive form
and moment matrix methods.

\begin{center}

\begin{algorithm}[\GIF -- \MMtx\;Method]
\label{alg:GIF-MMtx}

\noindent

\begin{itemize}
\item[]    \parskip0mm  \bf{Input:}  $ P = \{p_1,...,p_k\} \subseteq
\mathbb{R}[x_1,\ldots,x_n]$

\begin{itemize}
\item[] $Q_0 := P$
\item[] $ j:= 0 $
\item[]     \bf{do}
\textnormal{
\begin{itemize}
\item[]   $  d :=  \dim \ker \GIF (Q_j)  $
\item[]   $ Q_{j+1} := \gen( \GIF (Q_j) ) $
\item[]   $  r := \rank( \MMtx (Q_{j+1}) )$
\item[]   $ Q_{j+2} := \gen(\ker \MMtx (Q_{j+1}) ) $
\item[]   $j: = j + 2$
\end{itemize}}
\item[] \bf{until} $r = d $
\end{itemize}
\item[] \bf{Output}:  $ Q = \{q_1,...,q_\ell \} \subseteq
\mathbb{R}[x_1,\ldots,x_n]$  \\
  \hspace{1cm} $ Q $ \textnormal{ is in geometric involutive form } \\
  \hspace{1cm} $\sqrt[\mathbb{R}]{ \left\langle P
\right\rangle_\mathbb{R}  }  \; \;   \supseteq  \; \;  \left\langle Q
\right\rangle_\mathbb{R}  \; \;  \supseteq  \; \;  \left\langle P
\right\rangle_\mathbb{R}    $.
\end{itemize}

\end{algorithm}
\end{center}

\noindent
\textbf{Proof of the termination of Algorithm \ref{alg:GIF-MMtx}:}
We prove termination of the \GIF -- \MMtx\;Method under the assumption that suitable
generic points, if available, are determined at each iteration of the method.

\noindent
\textbf{Rank-Dim-Involutive Stopping Criterion}:
A natural termination criterion used in Algorithm \ref{alg:GIF-MMtx} is that the generators stabilize at some iteration and the system
is involutive:
\begin{equation}
\label{ConjMMtx}
\gen (\GIF (Q)) = \gen (\ker \MMtx (Q)) \; \mbox{and} \;  Q \; \mbox{involutive}
\end{equation}

Since different representations of the rings are involved we will focus on one, that of
polynomial generators during the proof.

In terms of generators our termination criterion $\rank (\MMtx (Q_{j+1})) = \dim \ker  \GIF( Q_j)$
 is expressed as $\gen (\GIF (Q_{j})) = \gen (\ker \MMtx (Q_{j+1}) )$.

Then  $\gen( \ker \MMtx(Q_{j+1}))$ and $\gen ( \GIF (Q_j))$ are both ideals of the system $P$.
 Since a generator of the geometric involutive form will also be a generator of the ideal in the moment matrix
at each iteration we have $\gen ( \GIF (Q_j)) \subseteq \gen ( \ker \MMtx(Q_{j+1}))$ in our algorithm.
Suppose the algorithm never stops,  then we will get a infinite ascending chain of ideals
with a strict inclusion at each iteration of the form $Q_j  \subset Q_{j+1}$ where
$ Q_j =\gen(\GIF (Q_{j-1}))$ and $Q_{j+1} =  \gen( \ker \MMtx(Q_j))$.
This is a violation of the ascending chain condition since $\mathbb{R}[x_1, ... , x_n]$ is a Noetherian Ring. Therefore, the generators must stabilize in the end and when stabilized, $Q$ is also involutive. 
 \hfill $\square$

The algorithm above uses the following subroutines.

\begin{center}

\begin{algorithm}[\GIF]
\label{alg:GIF-subroutine}

\noindent

\begin{itemize}
\item[]    \parskip0mm  \bf{Input:} $ Q \subseteq \mathbb{R}[x_1,\ldots,x_n]$
\item[] {\bf{Output}}:  \textnormal{Return a geometric involutive form $\GIF(Q)$.}
\end{itemize}

\end{algorithm}

\end{center}

Note the algorithm \ref{alg:ProjInvBasis} is an explicit implementation of \GIF. 

\begin{center}

\begin{algorithm}[\MMtx]
\label{alg:MMtx-subroutine}

\noindent

\begin{itemize}

\item[]    \parskip0mm  \bf{Input:} $ Q \subseteq \mathbb{R}[x_1,\ldots,x_n]$.  \textnormal{Set} $ d := \deg(Q)$.
            \textnormal{   \begin{enumerate}
                \item Construct the general $N(n,d) \times N(n,d)$  moment matrix.
			\item Construct the involutive prolongation $D^{d} Q$.
                \item Use SDP methods to numerically solve for a generic point that maximizes the rank of the moment matrix subject to the constraints $D^d Q$.
                \end{enumerate}}

\item[] {\bf Output}: \textnormal{ Return $\MMtx (Q) \succeq 0$ the moment matrix evaluated at this generic point.}
\end{itemize}

\end{algorithm}
\end{center}

\begin{center}

\begin{algorithm}[\gen]
\label{alg:gen-subroutine}

\noindent

\begin{itemize}
\item[]    \parskip0mm  \bf{Input:} \textnormal{
 $\GIF (Q)$ or $\ker \MMtx (Q) $}

\item[] {\bf Output}: \textnormal{  Polynomial generators corresponding to $\GIF (Q)$ or $\ker \MMtx (Q) $}
\end{itemize}

\end{algorithm}

\end{center}

\subsection{Two variable example}

\label{example:2vars}
Consider the polynomial system with two variables $x$ and $y$.
\begin{equation} \label{ex:2vars}
 P_2 =\{(y^2 -1)^2, (y^2 -1)(x^2 - 1) \}
\end{equation}

First we apply \GIF \;to $P_2$ to compute the involutive form of it. The dimension table is in Figure~\ref{TableCMtxSimpleEx2vars}.

\begin{figure}[h!]
\centerline{
\begin{tabular}{|l|cccc|}
\hline
            & $k=0$ & $k=1$ & $k=2$ & k=3    \\
\hline
$\ell=0$ &  13       &  15      &  17        &    19     \\
$\ell=1$ &  10    &   12  &  14        &    16    \\
$\ell=2$ &  6       &  9      &  \color{blue}{\boxed{11}} &    13    \\
$\ell=3$ &  3       &  6      &  9        &   11    \\
$\ell=4$ &  1       &  3      &  6        &    9    \\
\hline
\end{tabular}}
\caption{ Table of $\dim  \ProjKer^\ell  \ProKer^k (P_2)$ for system (\ref{ex:2vars})  The (blue) boxed $11$ in the third
		column corresponds to $\ProjKer^2 \ProKer^2 (P_2)$.
		\label{TableCMtxSimpleEx2vars}}
\end{figure}

Now $\dim \ProjKer^2 \ProKer^2 (P_2) = \dim  \ProjKer^3 \ProKer^3 (P_2) $ so $\ProjKer^2 \ProKer^2 (P_2) $ satisfies one of the conditions for
an involutive system.   The second condition is that the symbol of $\ProjKer^2 \ProKer^2 (P_2)$ is involutive.
Applying the symbol test~\eqref{projelimtest1} and we find that $\dim \Symbol \; \ProjKer^2 \ProKer^3 (P_2) = \dim \Symbol \; \ProjKer^2 \ProKer^2 (P_2 ) = 2 $,
so the symbol of it turns out to be involutive as well. Therefore $\ProjKer^2 \ProKer^2 (P_2)$ is involutive.

Now we apply the subroutine $\MMtx$\;to $\gen (\ProjKer^2 \ProKer^2 (P_2))$ to compute the moment matrix $M$. We convert $\ker M$ into polynomial generators by subroutine $\gen$.
The dimension of $\ker M$ is 6 which means there are 6 generators in $\gen(\ker M)$, which are moderately complicated
numerical polynomials.

We again apply $\GIF$ to $\gen(\ker M)$ to compute the involutive form. The dimension table is shown in Figure~\ref{TableCMtxSimpleEx2vars1}.
The input system corresponding to the (red) boxed 9 is already involutive. As mentioned in Remark \ref{rem:extracProj} in algorithm \ref{alg:ProjInvBasis} and more generally in \GIF\;algorithm, we can extract projected systems of lower degree than input system. This improves on our previous algorithm \cite{SRWZ2010}. We demonstrate this procedure here.   In Figure~\ref{TableCMtxSimpleEx2vars1}, the system corresponding to the red boxed 9 is involutive and has degree 4.
Since $N(2,4) = 15$  there are $15 -9 = 6$ polynomials in the system.  However descending further down the column of the table, we find the system corresponding to the blue boxed 5 is also involutive.   In that case $N(2,2) = 6$ so there is only $1$ corresponding generator.

\begin{figure}[h!]
\centerline{
\begin{tabular}{|l|cccc|}
\hline
            & $k=0$ & $k=1$ & $k=2$ & k=3    \\
\hline
$\ell=0$ &  \color{red}{\boxed{9}}       &  11      &  13       &    15     \\
$\ell=1$ &  7   &   9  &  11      &    13    \\
$\ell=2$ &   \color{blue}{\boxed{5}}     &  7    &  9  &   11    \\
$\ell=3$ &  3      &  5      &  7     &   9    \\
$\ell=4$ &  1      &  3      &  5     &   7    \\
\hline
\end{tabular}}
\caption{ Table of $\dim    \ProjKer^\ell  \ProKer^k \gen (\ker M)$ for the first \GIF -- \MMtx \;
iteration in Example \ref{example:2vars}.
 }
 \label{TableCMtxSimpleEx2vars1}
\end{figure}

We apply $\gen$ to compute the generator set:
\begin{equation}\label{gen:2vars}
  \{  0.7071067 * y^2 - 0.7071067 + \mbox{small terms less than}\; 10^{-11}   \}
\end{equation}

If we apply \GIF \;to equation \eqref{gen:2vars}, the dimension table is exactly the same as the one in Figure \ref{TableCMtxSimpleEx2vars1}.
Therefore the projected system is equivalent to the input system. After normalization and ignoring small terms, we get $y^2 - 1$ which is a geometric involutive basis for the real radical for $P_2$.

\subsection{Three variable example}

\label{sec:3varEx}

In this section we apply the \GIF -- \MMtx \; method to the following trivariate system with \GIF \;explicitly implemented by Algorithm \ref{alg:ProjInvBasis}:
\begin{equation}
\label{3varsys}
 P_3=\begin{cases}
               x^2y^2 - y^4 + y^2z^2 -x^2 -z^2 +1\\
                x^2y^2 - y^4 + y^2z^2 +x^2 -2y^2 + z^2 -1\\
               x^4z+x^2z^3 -2x^2y^2 -x^2z -z^3 -2x^2 +2y^2 +2 \\
               x^4z +x^2 z^3 -2 x^2y^2 +x^2 z + z^3 -2 x^2  -2 y^2 -2
            \end{cases}
\end{equation}

We first apply subroutine $\GIF$ to $P_3$. The dimension table is shown in Figure~\ref{TableCMtxSimpleEx3vars}.

\begin{figure}[h!]
\centerline{
\begin{tabular}{|l|ccccc|}
\hline
            & $k=0$ & $k=1$ & $k=2$ & $k=3$ & k=4  \\
\hline
$\ell=0$ &  46  &  57   &  66   &   73  &  79   \\
$\ell=1$ &  29  &  38   &  46   &   53  &  59   \\
$\ell=2$ &  \color{red}{\boxed{15}}  &  21   & \color{blue}{\boxed{27}}   &   33  &  39   \\
$\ell=3$ &  9   &  15   &  21   &   27  &  33   \\
$\ell=4$ &  4   &   9   &  15   &   21  &  27   \\
\hline
\end{tabular}}
\caption{ Table of $\dim     \ProjKer^\ell \ProKer^k (P_3)$ for system (\ref{3varsys}).
 }
 \label{TableCMtxSimpleEx3vars}
\end{figure}

At prolongation zero of Algorithm \ref{alg:ProjInvBasis} we determine if there are any projected involutive systems whose prolongations yield the same ideal as the system (so that the prolongations can be
discarded).
We find such an involutive system $\ProjKer^2 \ProKer^0 (P_3)$ which corresponds to the red boxed 15 in column 1 of Figure \ref{TableCMtxSimpleEx3vars}.
From the dimension information we can deduce that since the number of monomials of degree $\leq 3$
is $N(3, 3) = 20$ there will be $20 - 15 = 5$ polynomials generators corresponding to
$\ProjKer^2 \ProKer^0 (P_3)$.
System $\ProjKer^3 \ProKer^0 (P_3) = \bar{R}$ is of lower degree and also easily found to be involutive.  However it does not satisfy the inclusion test of  Algorithm \ref{alg:ProjInvBasis} given by
$\ProKer^{d + k-\bar{d}} \bar{R}  \subseteq \ProKer^k (P_3)$ which shows that it is not equivalent to the original system.
We find that $\ProjKer^2 \ProKer^0 (P_3)$ does satisfy the inclusion output condition, so we exit \GIF \;   and
apply subroutine $\MMtx$ to  $\gen(\ProjKer^2 \ProKer^0 (P_3))$.
In our previously published method we would have first identified
the blue boxed 27 corresponding to the involutive system $\ProjKer^2 \ProKer^2 (P_3)$.
Our approach is a clear improvement, and avoids creating the large degree 5 moment matrix of the previous
approach.

We compute the generator set of the moment matrix $M$ using the subroutine $\gen(\ker M)$.
The rank of moment matrix is $7$ which means $\gen(\ker M)$ has dimension $13$. We apply $\GIF$ \; to $\gen(\ker M)$ and the dimension table is given in Figure~\ref{TableCMtxSimpleEx3vars1}
\begin{figure}[h!]
\centerline{
\begin{tabular}{|l|cccc|}
\hline
            & $k=0$ & $k=1$ & $k=2$ & k=3   \\
\hline
$\ell=0$ &  7   &  9   &  11   &   13    \\
$\ell=1$ &  \color{blue}{\boxed{5}}   &  7   &  9   &   11    \\
$\ell=2$ &  3  &  5    &  7   &   9    \\
$\ell=3$ &  1   &  3    &  5    &   7    \\
\hline
\end{tabular}}
\caption{ Table of $\dim     \ProjKer^\ell  \ProKer^k \gen(\ker M)$ in the moment matrix
calculation for $\gen(\ProjKer^2 \ProKer^0 (P_3))$.
 }
 \label{TableCMtxSimpleEx3vars1}
\end{figure}

In this iteration of $\GIF$ three systems are involutive and correspond to the
$\ell = 0, 1, 2$ entries of column 1 of  Figure \ref{TableCMtxSimpleEx3vars1}.
Corresponding to the elimination of higher order systems by the inclusion test in Algorithm \ref{alg:ProjInvBasis}, 
we can discard 2 of the 3 systems, which correspond to $\ell = 0$ and $\ell = 2$ entries in the first column.
The output lower degree geometric involutive basis therefore corresponds to the blue boxed entry
in the figure.

At the next iteration the generators corresponding to $\ell = 1$ are sent to the moment matrix.
We find that the termination condition is satisfied, that is $d = 5 = r$.
The algorithm then terminates with an output of $10 - 5 = 5 $ generators.

To get more insight into the output we now analyze it further.
From the Figure \ref{TableCMtxSimpleEx3vars1} we see that there is a projected system corresponding to
$\ell = 2$ of dimension $3$ and degree $1$.
When it is extracted we find a single nice generator:
\begin{equation} \label{ex:3var1}
0.8944271-0.4472135z+5.5511151 \times 10^{-17} * y
\end{equation}
After dropping off the small term and normalization, we get $z - 2$.
Now we consider the other generators of degree $2$.
We can simplify them by substituting $z = 2$ from the projected generator and find
\begin{equation} \label{ex:3vars2}
-0.3015113*x^2+0.3015113*y^2-0.9045340 + \mbox{small terms less than}\; 10^{-15}.
\end{equation}
which is approximately $x^2 -y^2 +3$.
Thus our output geometric involutive basis is
\begin{equation}
\label{ex3var:outputGIFbasis}
 \{z - 2, x(z-2), y(z-2), z(z-2),x^2 - y^2 + 3\}
\end{equation}
A hand calculation checks that this is a geometric involutive basis for the real radical of the input system.

\section{Discussion}

In this paper we present improvements of our numerical geometric
involutive bases for polynomial systems of equations.
We also began an exploration of the interaction of these methods with
SPD programming methods and computation of such bases for positive dimensional real radical ideals.

We give methods to extract and decrease the degree of
immediate systems and the output basis.  One such tool is an inclusion test
whereby higher degree redundant systems can be discarded.
Prompted by a number of requests we have given more details of our implementation
of Cartan's involutivity  test for positive dimensional ideals.
Reduction of degree techniques are critical and have been extensively developed in the
symbolic case for Gr\"obner bases \cite{Faugere99} and triangular decompositions
\cite{Chen13,CMXY09}.
Significant progress has also been made in symbolic-numeric methods such
as border bases \cite{Mourrain:1996,Mourrain:1999,MourrainTrebuchet:2000,Mourrain:2002}
in removing higher degree polynomials. 
Perhaps the closest objects to geometric involutive bases in the zero dimensional case
are H-Bases \cite{MollerSauer:2000}.

Moreover, we were motivated by remarkable recent work by Lasserre and
collaborators \cite{SemiDefinite:2012} using SDP methods
for identifying the real radical of zero dimensional polynomial ideals.

The work  \cite{SemiDefinite:2012} motivated us to combine
SDP -- moment matrix methods with our geometric involutive bases to approximate
positive dimensional real radical ideals.
In particular,  the termination criterion $\rank(\MMtx(Q)) = \dim \ker \GIF(Q) $ in Algorithm \ref{alg:ProjInvBasis}  is equivalent to the rank stabilization condition  in Lasserre \cite{SemiDefinite:2012} for zero dimensional systems.  Moreover in our initial explorative experiments we obtained generators for the real radical of positive dimensional ideals for a small set of examples and deserves further study.

In our preliminary study in order to study the interaction between these two
methods we focused on an algorithm that cleanly separates the step of taking
a geometric involutive basis at each iteration of the algorithm.
An alterative strategy that we will pursue in future work is motivated
by the approach
of Lasserre et al in the zero dimensional case \cite{LasserreLaurentRostalski09}.  Instead of demanding
a (projected) involutive
form at each iteration, they allowed the iteration and prolongation of
moment matrices until
the projected criteria for involution were obtained (that is a zero
dimensional symbol in that case).
This has the advantage that geometric involutive form calculations
whose complexity implicitly depends
on the total number of complex solutions are avoided until later, when
such complex solutions have
been discarded as a result of new generators being found in the kernel
of the moment matrix.


\begin{thebibliography}{999}

\small{


\bibitem{SemiDefinite:2012}
Anjos, Miguel F., Lasserre, Jean B., 2011. Handbook on Semidefinite,
Conic and Polynomial Optimization. Springer, 25--49.

\bibitem{BasuPollackRoy06} S. Basu, R. Pollack and M-F Roy,
~Algorithms in real algebraic geometry. 2nd edition. Algorithms and
Computation in Math., 10. Springer-Verlag, 2006.

\bibitem{BHSW13} D.J. Bates, J.D. Hauenstein, A.J. Sommese, and C.W. Wampler,
Numerically Solving Polynomial Systems with the Software Package
Bertini. In preparation. To be published by SIAM, 2013.

\bibitem{BDHSW12} G.M Besana, S. DiRocco, J.D. Hauenstein, A.J.
Sommese, and C.W. Wampler. Cell decomposition of almost smooth real
algebraic surfaces. Num. Algorithms, DOI:10.1007/s11075-012-9646-y,
(published online Sept 28, 2012).

\bibitem{BPT2012}
Semidefinite Optimization and Convex Algebraic Geometry,
Editors: Grigoriy Blekherman, Pablo A. Parrilo and Rekha R. Thomas,
MOS-SIAM Series on Optimization 13 , 2012.


\bibitem{BLRSZ:2004}
Bonasia, J., Lemaire, F., Reid, G., Scott, R., Zhi, L.. Determination of
approximate symmetries of differential equations. In: CRM Proceedings and
Lecture Notes. pp. 233--249. 2004.


\bibitem{Chen13} C. Chen, J.H. Davenport, J.P. May, M.M. Maza, B. Xia and
R. Xiao.
\newblock {\em ~Triangular decomposition of semi-algebraic systems}
\newblock Journal of Symbolic Computation Vol49-0, pp 3-26, 2013.

\bibitem{CMXY09} C. Chen, M. Moreno Maza, B. Xia, L. Yang, ~Computing
Cylindrical Algebraic Decomposition via Triangular Decomposition.
Proc. of ISSAC 2009, pages 95-102, ACM Press, New York, 2009.

\bibitem{Collins75} G. Collins, Quantifier elimination for real closed
fields by cylindrical algebraic decomposition. Springer Lec. Notes
Comp. Sci. v33. 515-532, 1975.

\bibitem{CurtoFialkow96} R. Curto, L. Fialkow. Solution of the
truncated complex moment problem for
at data. Mem. Amer. Math. Soc. 119 (1996), no. 568.

\bibitem{DavenportHeintz88} J.H. Davenport and J. Heintz,
Real quantifier elimination is doubly exponential,
J. Symbolic Comp., 5:29-35, 1988.

\bibitem{Faugere99}
J.-C. Faug\`{e}re. A new efficient algorithm for computing Gr\"obner
bases (F4). Journal of Pure and Applied Algebra, 139(1-3):61-88, 1999.

\bibitem{GerdtBlinkov:1998}
Gerdt, V., Blinkov, Y., 1998. Involutive bases of polynomial ideals.
Mathematics and Computers in Simulation 45, 519--541.

\bibitem{GTZ88} Gianni, P., Trager, B., Zacharias, G., 1988.
Gr\"obner Bases and. Primary Decomposition of Polynomial Ideals. J.
Symb. Comp. 6, 149-167.

\bibitem{GrigorevVorobjov92} D. Grigor\'{e}v and N. Vorobjov.
~Counting connected components 
subexponential time. Comput. Complexity, 2(2), 133-186, 1992.


\bibitem{Hauenstein12} J. Hauenstein, Numerically computing real
points on algebraic sets. Acta Appl. Math.
DOI:10.1007/s10440-012-9782-3, (published online September 27, 2012).

\bibitem{Hong90} H. Hong. ~Improvement in CAD-Based Quantifer Elimination.
Ph.D. thesis, the Ohio State University, Columbus, Ohio, 1990.

\bibitem{Kuranishi:1957}
Kuranishi, M., 1957. On {E}. {C}artan's prolongation theorem of exterior
differential systems. Amer. J. Math 79, 1--47.

\bibitem{LasserreLaurentRostalski09} J.B. Lasserre, M. Laurent and P. Rostalski.
~A prolongation-projection algorithm for computing the finite real
variety of an ideal. Theoret. Comput. Sci., 410(27-29), 2685-2700,
2009.

\bibitem{Lu06} Y. Lu, Finding all real solutions of polynomial
systems, Ph.D Thesis, University of Notre Dame, 2006.
Results of this thesis appear in:\\
{\small (with D.J. Bates, A.J. Sommese, and C.W. Wampler),
Finding all real points of a complex curve, Contemporary
Mathematics 448 (2007), 183--205.}

\bibitem{Macaulay:1916}
Macaulay, F., 1916. The Algebraic Theory of Modular Systems. Vol.~19. Cambridge
tracts in Math. and Math. Physics.

\bibitem{MaZhi12} Y. Ma and L. Zhi.
~Computing Real Solutions of Polynomial Systems via Low-rank Moment
Matrix Completion
Proc. 2012 Internat. Symp. Symbolic Algebraic Comput. pp. 249-256.


\bibitem{MollerSauer:2000}
M{\"o}ller, H., Sauer, T., 2000. H-bases for polynomial interpolation and
system solving. Advances Comput. Math. 12, 23--35.

\bibitem{Mourrain:1996}
Mourrain, B., 1996. Isolated points, duality and residues. J. of Pure and
Applied Algebra 117 {\&} 118, 469--493.

\bibitem{Mourrain:1999}
Mourrain, B., 1999. A new criterion for normal form algorithms. In: Fossorier,
M., Imai, H., Lin, S., Poli, A. (Eds.), AAECC. Vol. 1719. Springer, Berlin,
pp. 430--443.


\bibitem{MourrainTrebuchet:2000}
Mourrain, B., Tr{\'e}buchet, P., 2000. Solving projective complete intersection
faster. In: Traverso, C. (Ed.), Proc. 2000 Internat. Symp. Symbolic Algebraic
Comput. {ISSAC}'00. ACM Press, New York, pp. 430--443.

\bibitem{Mourrain:2002}
Mourrain, B., Tr{\'e}buchet, P., 2002. Algebraic methods for numerical solving.
In: Proceedings of the 3rd International Workshop on Symbolic and Numeric
Algorithms for Scientific Computing. pp. 42--47.

\bibitem{Pommaret:1978}
Pommaret, J., 1978. Systems of Partial Differential Equations and Lie
Pseudogroups. Gordon and Breach Science Publishers.

\bibitem{RRS00} F. Rouillier, M.-F. Roy, and M. Safey El Din. ~Finding
at least one point in each connected component of a real algebraic set
defined by a single equation. J. Complexity, 16 (4), 716-750, 2000.

\bibitem{SRWZ2010}
Scott, R., Reid, G., Wu, W., Zhi, L.
Geometric Involutive Bases and Applications to Approximate Commutative Algebra.
In Approximate Commutative Algebra, Eds Robbiano, L. and Abbott, J.
Texts and Monographs in Symbolic Computation, pp. 99 --124, 2010.

\bibitem{ReidLinWittkopf:2001}
Reid, G., Lin, P., Wittkopf, A., 2001. Differential elimination-completion
algorithms for {DAE} and {PDAE}. Studies in Applied Mathematics 106~(1),
1--45.

\bibitem{ReidWittkopf:2001}
Reid, G., Wittkopf, A., 2001.
Fast Differential Elimination in {C}: The {CD}iffElim {E}nvironment.
Comp. Phys. Comm. 139~(2), 192--217.

\bibitem{ReidSmithVerschelde:2002}
Reid, G., Smith, C., Verschelde, J., 2002. Geometric completion of differential
systems using numeric-symbolic continuation. SIGSAM Bulletin 36~(2), 1--17.

\bibitem{ReidTangZhi:2003}
Reid, G., Tang, J., Zhi, L., 2003. A complete symbolic-numeric linear method
for camera pose determination. In: Sendra, J. (Ed.), Proc. 2003 Internat.
Symp. Symbolic Algebraic Comput. {ISSAC}'03. ACM Press, New York, pp.
215--223.

\bibitem{ReidZhi2009}
G. Reid and L. Zhi,  Solving polynomial systems via symbolic-numeric reduction to geometric involutive form.
J. Symb. Comput.  pp 280--291.  2009.

\bibitem{Seidenberg54} A. Seidenberg, ~A new decision method for
elementary algebra, Ann. of Math. 60, 365-374, 1954.

\bibitem{Seiler2010} W. M. Seiler.
\newblock {\em ~Involution. The formal theory of differential
equations and its applications in computer algebra}.
\newblock Springer Verlag, Algorithms and Computation in Mathematics.
Vol 24. 2010.


\bibitem{SommeseWampler05} A.J. Sommese and C.W. Wampler.
\newblock {\em ~The Numerical solution of systems of
polynomials arising in engineering and science}.
\newblock World Scientific Press, 2005.

\bibitem{Stetter:2004}
Stetter, H., 2004. Numerical Polynomial Algebra. Soc. for Industrial
and Applied Math. (SIAM).

\bibitem{Tarski31} A. Tarski, ~A decision method for elementary
algebra and geometry,
Fund. Math. 17, 210-39, 1931.









\bibitem{WuZhi12} X. Wu and L. Zhi.
~Computing the multiplicity structure from geometric involutive form,
Journal of Symbolic Computation, 47(3): 227-238, 2012.

\bibitem{WuReid13} Wenyuan Wu and Greg Reid.
~Finding points on real solution components and applications to
differential polynomial systems. ISSAC 2013: 339-346.

\bibitem{YZZ13}
Zhengfeng Yang, Lihong Zhi, and Yijun Zhu.
~Verfied error bounds for real solutions of positive-dimensional
polynomial systems
In ISSAC'2013 Proc. 2013 Internat. Symp. Symbolic Algebraic Comput. pp. 371-378.
}

\end{thebibliography}
\end{document}